\documentclass[11pt]{article}
\begin{document}
\title{On the {\'e}tale cohomology of algebraic varieties with totally
degenerate reduction over $p$-adic fields\\\textit{{\small to J.
Tate}}}
\author{Wayne Raskind\thanks{Partially supported by NSF grant 0070850,
SFB 478 (M{\"u}nster), CNRS France, and sabbatical leave from the
University of Southern California}\hspace{.1in}
\setcounter{footnote}{6} Xavier Xarles\thanks{Partially supported by
grant BHA2000-0180 from DGI}}

\date{}
\maketitle

\newcommand{\free}[1]{{#1}/tors}

\newcommand{\bC}{{\bf C}}
\newcommand{\bZ}{{\bf Z}}
\newcommand{\bQl}{{\bf Q}_{\ell}}
\newcommand{\bZl}{{\bf Z}_{\ell}}
\newcommand{\Xb}{{\overline X}}
\newcommand{\Zb}{{\overline Z}}
\newcommand{\Eb}{{\overline E}}
\newcommand{\Kb}{{\overline K}}
\newcommand{\Cb}{{\overline C}}

\newcommand{\cX}{{\cal X}}
\newcommand{\cXb}{{\overline{\cal X}}}
\newcommand{\Yb}{{\overline Y}}
\newcommand{\bQ}{{\bf Q}}
\newcommand{\bQp}{{\bf Q}_p}
\newcommand{\bZp}{{\bf Z}_p}
\newcommand{\bCp}{{\bf C}_p}
\newcommand{\bR}{{\bf R}}
\newcommand{\cE}{{\cal E}}
\newtheorem{proposition}{Proposition}
\newtheorem{theorem}{Theorem}
\newtheorem{lemma}{Lemma}
\newtheorem{definition}{Definition}
\newtheorem{remark}{Remark}
\newtheorem{example}{Example}
\newtheorem{corollary}{Corollary}
\newtheorem{conjecture}{Conjecture}
\newtheorem{observation}{Observation}
\newcommand{\plim}{\displaystyle{\lim_{\stackrel{\longleftarrow}{n}}}\,}
\newcommand{\bF}{{\overline F}}
\newcommand{\kb}{{\overline k}}
\newcommand{\Fb}{{\overline F}}
\newcommand{\cO}{{\cal O}}
\newcommand{\im}{\mbox{Im}}
\newcommand{\bG}{{\bf G}}
\begin{abstract}  Let  $K$ be a field of characteristic zero that is complete with respect to
a discrete valuation, and with perfect residue field.  We formulate
the notion of {\it totally degenerate reduction} for a smooth
projective variety $X$ over $K$.
 We show that
for all prime numbers $\ell$, the $\bQl$-\'etale cohomology of such
a variety is (after passing to a finite unramified extension of $K$)
a successive extension of direct sums of Galois modules of the form
$\bQl(r)$.  More precisely, this cohomology has an increasing
filtration whose $r$-th graded quotient is of the form
$V\otimes_{\bQ}\bQl(r)$, where $V$ is a finite dimensional
$\bQ$-vector space that is independent of $\ell$, with an unramified
action of the absolute Galois group of $K$.
\end{abstract}
\parindent=0cm
{\bf AMS Subject Classification}:  14F20, 14F30 (primary) and 14G20
(secondary)
\section*{Introduction}

Let $K$ be a field of characteristic zero that is complete with
respect to a discrete valuation, and let $A$ be an abelian variety
of dimension $d$ over $K$.  Let $\Kb$ be an algebraic closure of
$K$, $\ell$ be a prime number and $A[\ell^n]$ be the group of points
of $A$ of order $\ell^n$ over $\Kb$.  We set

$$T_{\ell}(A)=\plim A[\ell^n].$$

Assume that $A$ has totally split multiplicative reduction, or
equivalently, that the connected component of the identity of the
special fibre of the N\'eron model of $A$ over the ring of integers
of $K$ is a split torus. Then the theory of $p$-adic uniformization
of Tate (unpublished), Mumford [Mu2] and Raynaud [Ray1] implies that
$A$ can be realized as the rigid analytic quotient of a
multiplicative torus by a lattice $\Lambda$. This implies easily
that for {\it every} $\ell$, there is an extension of Galois
modules:

$$0\to \Lambda'\otimes_{\bZ}\bZl(1)\to T_{\ell}(A)\to \Lambda\otimes_{\bZ}\bZl\to 0,$$

where $\Lambda'$ is a lattice that is canonically isogenous to
$\Lambda$.  In this paper, we generalize these results to describe
in such terms the \'etale cohomology of a class of varieties having
a type of reduction that we call {\it totally degenerate}. We employ
purely algebro-geometric methods, mainly usage of the comparison
theorem between $p$-adic \'etale cohomology and log-crystalline
cohomology (the semi-stable conjecture $C_{st}$ proven by Tsuji
[Tsu]), to prove our results. While many examples of these varieties
have uniformizations, one does not expect that all such do, and we
feel that our methods will complement uniformization methods even
for the case of abelian
varieties as above.\\

This is the first in a series of papers where we study closely the
cohomology of varieties with totally degenerate reduction.  In [RX],
we apply the results of this paper to define and study ``$p$-adic
intermediate Jacobians.''  In [R2], we formulate a conjecture of
Hodge-Tate type for such varieties, which would describe the
coniveau filtration on $p$-adic cohomology in terms of the kernels
of ``enriched monodromy operators.'' In [R1], we prove a form of
this conjecture for divisors.  The second author and Infante are
also studying complex analogues of some of the
results in this paper.\\

To describe more precisely the contents of this paper, let $X$ be a
smooth, projective, geometrically connected variety over $K$. Let
$\Xb$ denote base extension of $X$ to $\Kb$. We begin in \S 2 by
formulating the notion of totally degenerate reduction.  While this
term is used often, we have not seen a precise general definition,
and formulating such is one important part of this paper.  We view
it as a set of conditions that are to be satisfied by the components
and their various intersections of the special fibre of a suitable
regular proper model of $X$ over the valuation ring $R$ of $K$ with
strictly semi-stable reduction (we assume that such a model exists).
Roughly speaking, the conditions say that these intersections are
very simple cohomologically.  Our main result (see Corollary 1 of \S
4 and Theorem 3 of \S 6) says that for {\it all} $\ell$, the \'etale
cohomology groups $H^*(\Xb,\bQl)$ are (after passing to a finite
unramified extension of $K$) successive extensions of $\bQl$ by
$\bQl(r)$ for suitable $r$. More precisely, we show that there is a
good {\it monodromy filtration}, whose graded quotients have a
$\bZ$-structure that is given in a natural way by the cohomology of
the {\it Chow complex} that is formed from the Chow groups
 of the components of the special fibre and their
intersections (see \S 3 for the definition of this complex).
 \\

For $\ell\neq p$, our main result is essentially known, as it
follows without too much difficulty from the work of Rapoport-Zink
[RZ1], although the result is not explicitly stated there (see \S 4
for more details). It is the $p$-adic cohomology part that requires
a careful analysis of several different filtrations and their
eventual coincidence. To do this, we discuss  in \S 5 the monodromy
filtration defined by Mokrane [Mok] on the log-crystalline
cohomology of the special fibre of a regular proper model of $X$
over $R$ with semi-stable reduction, as defined by Hyodo-Kato [HK].
Then, in \S 6, we use work of Hyodo [H] and Tsuji [Tsu] to ``lift''
this filtration to a monodromy filtration on the $p$-adic {\'e}tale
cohomology of $\Xb$. We note that there is no ``simple'' monodromy
filtration on $p$-adic cohomology, in general, as was first pointed
out by Jannsen [Ja1], and we can only expect to get such a
filtration on $\ell$-adic cohomology for all $\ell$ for some classes
of varieties like those considered here.  We feel that the existence
of such a filtration on $p$-adic cohomology is a very important and
useful result, which is the foundation for all
of our work in this direction.\\

  Examples of varieties to which the methods of this paper may be applied
include abelian varieties with totally multiplicative reduction and
products of Mumford curves [Mu1] or other $p$-adically uniformizable
varieties, such as Drinfeld modular varieties [Mus] and some unitary
Shimura varieties (see e.g. [La], [Z], [RZ2] and [Va]).
Unfortunately, these assumptions only include rather simple
varieties $X$ with good reduction, as in that case the \'etale
cohomology groups $H^*(\Xb,\bQl)$ are pure Galois modules for
$\ell\neq p$ by Deligne's theorem (Riemann hypothesis [De2]). Also,
there are some inequalities for the Hodge numbers that must be
satisfied for a variety with totally degenerate reduction and this
excludes certain types of varieties such as rigid Calabi-Yau 3-folds
(see Remark 6(ii) below for more details).  As pointed out to us by
Fontaine, there should be many more examples if one considers {\it
motives} with
totally degenerate reduction.  \\

        The authors would like to thank, respectively, the Universitat
Aut{\`o}noma
 de Barcelona and the University of Southern California
for their hospitality.  This paper was completed while the first
author enjoyed the hospitality of Universit{\'e} de Paris-Sud, SFB
Heidelberg and Universit{\'e} Louis Pasteur (Strasbourg). We also
thank J.-L. Colliot-Th{\'e}l{\`e}ne, L. Clozel, L. Fargues, J.-M.
Fontaine, O. Gabber, L. Illusie, K. K{\"u}nnemann, G. Laumon, A.
Quir\'os and T. Tsuji for helpful comments and information. Finally,
our hearty thanks to the referee, whose detailed comments have
greatly improved the contents and exposition of this paper.

 \section{Notation and Preliminaries}
         Let $K$ be a field of characteristic zero and complete with
respect to a discrete valuation, with valuation ring $R$ and perfect
residue field $F$ of characteristic $p>0$. We denote by $\Kb$ an
algebraic closure of $K$, by $\bF$ the residue field of $\Kb$, which
is an algebraic closure of $F$, by $W(F)$ the ring of Witt vectors
of $F$, by $W=W(\bF)$ the ring of Witt vectors of $\bF$, and by
$W_n=W_n(\bF)$ the ring of Witt vectors of length $n$. Denote by
$K_0$ the fraction field of $W(F)$ and by $L$ the fraction field of
$W$. The absolute Galois group $Gal(\Kb/K)$ will be denoted by $G$.
Note that we have a natural epimorphism $G \to Gal(\bF/F)$, so any
$Gal(\bF/F)$-module is naturally a $G$-module.
\\

For $X$ a smooth projective variety over $K$ we denote by $\Xb$ the
variety over $\Kb$ given by $X\times_K\Kb$, and in the same way for
a scheme $Z$ over $F$ we denote by $\Zb$ the scheme $Z\times_F\Fb$.
For $\ell$ a prime number and $r$ a nonnegative integer, we denote
by $\bQl(r)$ the Galois module $\bQl$, twisted $r$ times by the
cyclotomic
character.  If $r$ is negative, then $\bQl(r)=\mbox{Hom}(\bQl(-r),\bQl)$.\\

We denote by $B_{DR}$ the ring of $p$-adic periods of Fontaine. Fix
a uniformizer $\pi$ of $K$ and the extension of the $p$-adic
logarithm to $\Kb^*$ with $\log(\pi)=0$. These choices determine an
embedding of $B_{st}$, the ring of periods of semi-stable varieties,
in $B_{DR}$. See ([I1], \S 1.2) for definitions of these
rings and more details.\\

 Let $S$ be any domain with field of fractions $\mbox{frac}(S)$ of
characteristic zero. If $M$ is an $S$-module, let $M_{tors}$ be the
torsion subgroup of $M$. We denote by $\free{M}$ the torsion free
quotient of $M$; that is, $M/M_{tors}$. We say that a map between
$S$-modules is an isomorphism modulo torsion if it induces an
isomorphism between the torsion free quotients. If $M$ and $M'$ are
torsion free $S$-modules of finite rank, we say that a map
$\phi\colon M\to M'$ is an isogeny if it is injective with cokernel
of finite exponent (as abelian group). In this case, there exists a
unique map $\psi\colon M'\to M$, the dual isogeny, such that
$\phi\circ \psi=[e]$ and $\psi\circ \phi=[e]$, where $e$ is the
exponent of the cokernel of $\phi$ and $[e]$ denotes the map
multiplication by $e$. In general, if $M$ and $M'$ are $S$-modules
with torsion free quotients of finite rank, we say that a morphism
$\psi\colon M'\to M$ is an isogeny if the induced map on the torsion
free quotients is
an isogeny.\\

We will use subindices for increasing filtrations and superindices
for decreasing filtrations. If $M_{\bullet}$ is an increasing
filtration on an abelian group $H$ (respectively $F^{\bullet}$ a
decreasing filtration on $H$), we will denote by $Gr^M_{i}(H)$ the
quotient $M_{i}(H)/M_{i-1}(H)$ (respectively, $Gr_F^i(H)$ the
quotient $F^i(H)/F^{i+1}(H)$). Observe that, if $M_{\bullet}$ is an
increasing filtration, then
$F^{i}(H):=M_{-i}(H)$ is a decreasing filtration.\\

Let $X$ be a smooth projective geometrically connected variety over
$K$ of dimension $d$. We assume that $X$ has a regular proper model
$\cX$ over $R$ which is strictly
 semi-stable, which means that the following conditions hold:\\

(*)   Let $Y$ be the special fibre of $\cX$.  Then $Y$ is reduced;
write
 $$Y=\bigcup_{i=1}^{n}Y_i,$$
 with each $Y_i$ irreducible. For each nonempty subset
$I=\{i_1,\dots,i_k\}$ of $\{1,\dots,n\}$, we
 set
 $$Y_I=Y_{i_1}\cap...\cap Y_{i_k},$$
scheme theoretically.  Then $Y_I$ is smooth over $F$ of pure
codimension $|I|$ in $\cX$ if it is nonempty. See ([dJ] 2.16 and
[Ku1], \S 1.9, 1.10) for a clear summary of these conditions, as
well as comparison with other notions of semi-stability.\\

\section{Totally degenerate varieties and totally degenerate
reduction}

\begin{definition}
With notation as above, we say that $Y$ is {\it totally degenerate
over $F$} if there exists an embedding of $Y$ into a projective
space such that the following conditions are satisfied for each
$Y_I$ (where $d_I=\dim(Y_I)$):

 \begin{itemize}
 \item[a)] For every $i=0,...,d_I,$ the Chow groups $CH^i(\Yb_I)$ are
 finitely generated abelian groups.

The groups $CH^i_{\bQ}(\Yb_I)$ satisfy the Hodge index
 theorem:  let
 $$\xi:  CH^i(\Yb_I)\to CH^{i+1}(\Yb_I)$$ be the map given by
 intersecting with the class of a
 hyperplane section (which is determined by the fixed embedding of $Y$ into a projective space).
For every $i\le \frac {d_I}{2}$ and $x\in CH^i_{\bQ}(\Yb_I)$ such
that $\xi ^{d_I-2i+1}(x)=0$, we have that
$(-1)^i\mbox{deg}(x\xi^{d_{I}-2i}(x))\geq 0$, and equality holds if
and only if $x=0$. Here $\mbox{deg}$ denotes the degree map
$$\mbox{deg}: CH^{d_I}_{\bQ}(\Yb_I)\to \bQ$$

 \item[b)] For every prime number $\ell$ different from $p$, the
 {\'e}tale cohomology groups $H^{2i+1}(\Yb_I, \bZl)$ are torsion, and the
cycle map
 induces an isomorphism
 $$\free{CH^i(\Yb_I)}\otimes \bZl \cong \free{H^{2i}(\Yb_I, \bZl (i))}.$$

 Note that this is compatible with the action of the
 absolute Galois group $Gal(\Fb/F)$.

 \item[c)]  Denote by $H^*_{\rm crys}(\Yb_I/W)$ the crystalline cohomology groups of $\Yb_I$. Then
 $H^{2i+1}_{\rm crys}(\Yb_I/W)$ are torsion, and
 $$\free{CH^i(\Yb_I)\otimes W(-i) }\cong \free{H^{2i}_{\rm
 crys}(\Yb_I/W)}$$
 via the cycle map.  Here $W(-i)$ is $W$ with the action of Frobenius
multiplied by $p^{i}$.

 \item[d)]  $Y$ is ordinary, in that $H^r(\Yb,B\omega^s)=0$ for all $r$ and
 $s$. Here $B\omega$ is the subcomplex of exact forms in the logarithmic de
Rham
 complex on $\Yb$ (see e.g. [I3], D{\'e}finition 1.4).  By [H]
 and ([I2], Proposition 1.10), this is implied by the $\Yb_I$ being ordinary in the
 usual sense, in that $H^r(\Yb_I,d\Omega^s)=0$ for all $I, r$ and
 $s$. For more on the condition of ordinary, see ([I2],
 Appendice and [BK], Proposition 7.3)
\end{itemize}
\end{definition}

If the natural maps

$$CH^i(Y_I)\to CH^i(\Yb_I)$$

are isomorphisms modulo torsion for all $I$, we shall say that $Y$
is \textit{split totally degenerate}. Since the $CH^i(\Yb_I)$ are
all finitely generated abelian groups and there is a finite number
of them, there is a finite extension of the field of $F$ where all
the cycle classes in $CH^i(\Yb_I)$ modulo torsion are defined. So,
after a finite extension, any totally degenerate variety becomes
split totally degenerate.

\begin{remark}
\begin{itemize}
\item[(i)] For totally degenerate $Y$, the Chow groups
$CH^i(\Yb_I)_{\bQ}:=CH^i(\Yb_I)\otimes \bQ$
  satisfy the hard Lefschetz Theorem. That is, if $L$ is the class of a
  hyperplane section
  in $CH^1_{\bQ}(\Yb_I)$ considered in the Definition 1 a)
  and $\xi:CH^i_{\bQ}(\Yb_I)\to CH^{i+1}_{\bQ}(\Yb_I)$
  denotes the Lefschetz operator associated to $L$, then
  $\xi^{d_I-2i}:CH^i_{\bQ}(\Yb_I)\to CH^{d_I-i}_{\bQ}(\Yb_I)$
  is an isomorphism for all $i\le \frac {d_I}2$.  This follows from the Hard Lefschetz
Theorem in $\ell$-adic {\'e}tale cohomology, as proved by Deligne
[De3] and the bijectivity of the cycle map modulo torsion.

\item[(ii)]  Note that condition c) implies that for $Y_I$, we
have an isomorphism as $K_0$-vector spaces

$$(CH^i(\Yb_I)\otimes_{\bZ}L)^{Gal(\bF/F)}\to
H^{2i}(Y_I/W(F))\otimes_{W(F)}K_0(i)$$

where $K_0$ and $L$ are the fraction fields of $W(F)$ and $W$
respectively, and we are considering the diagonal action by
$Gal(\bF/F)$ on $CH^i(\Yb_I)\otimes_{\bZ}L$.

Using the well-known equivalence of categories between the category
of $p$-adic representations $V$ of $Gal(\bF/F)$ and the category of
finite dimensional $K_0$-vector spaces $D$ with a semilinear
endomorphism $\varphi$ whose slopes are all zero, the above
isomorphism implies that the $p$-adic representation
$CH^i(\Yb_I)\otimes_{\bZ}\bQp$ of $Gal(\bF/F)$ corresponds to the
$\varphi$-module $H^{2i}(Y_I/W(F))\otimes_{W(F)}K_0(i)$ via this
equivalence of categories.
\end{itemize}
\end{remark}

We will say that $X$ has \textit{totally degenerate reduction} if it
has a regular proper model $\cX$ over $R$ which is strictly
semi-stable and whose special fibre $Y$ is totally degenerate over
$F$. We find the name {\it totally degenerate reduction} a bit
pejorative, because one of the main themes of this paper is that
``bad reduction is good.''  However, since this terminology is
well-established, at least for dimension one, we decided to continue
with it.

\begin{remark}
The condition that the cycle map be an isomorphism modulo torsion
can be weakened to say that this map is an isomorphism when tensored
with $\bQl$ for all $\ell$, and the same with the crystalline cycle
map.  The results of this paper would be weakened to only give
$\bZ$-structures to the graded quotients up to isogeny, which are in
fact isomorphisms for almost all $\ell$ using the following lemma.
\end{remark}

 \begin{lemma}
   Let $Z$ be a smooth, projective, irreducible variety of dimension $d$ over a separably
closed field.  Assume that the Chow groups are finitely generated
abelian groups and that
   the cycle map is an isomorphism when tensored with $\bQl$ for all
$\ell$ different from the characteristic of the field.  Then for
almost all $\ell$, the integral cycle map:
 $$c_{\ell}:\: CH^i(Z)\otimes_{\bZ}{\bZl}\to H^{2i}(Z,\bZl(i))$$
 is an isomorphism.
 \end{lemma}
  {\bf Proof:} From the conditions, we know that the kernel and cokernel of the cycle map
 are finite groups.  Since the Chow groups are assumed to be finitely
generated abelian groups, the $\ell$-torsion is zero for almost all
$\ell$, and hence the map is injective for almost all $\ell$. As for
the cokernel, by a theorem of Gabber [G], the {\'e}tale cohomology
 groups $H^j(Z,\bZl)$ are torsion free for almost all $\ell$.  Consider
the following commutative diagram of pairings of finitely generated
$\bZl$-modules:

 $$\matrix{CH^i(Z)\otimes\bZl & \times & CH^{d-i}(Z)\otimes\bZl&\to& \bZ\cr
\downarrow&&\downarrow&&\downarrow\cr
 H^{2i}(Z,\bZl(i))&\times &H^{2d-2i}(Z,\bZl(d-i))&\to& \bZ_{\ell}}.$$
 Here the top row is the intersection pairing, tensored with $\bZl$, and the bottom pairing is
 cup-product on cohomology.  The vertical maps are the cycle maps.
By our assumptions on $Z$ and Poincar{\'e} duality, the
 intersection pairing on the (integral) Chow groups is perfect when
 tensored with $\bQ$, and so its determinant is a
 non-zero integer, say $m$. Let $S$ be the finite set of prime numbers consisting of those $\ell$ such that $H^{2i}(Z,\bZl)$
 or $H^{2d-2i}(Z,\bZl)$ has
 torsion, $\ell$ divides $m$, or $CH^i(Z)$ or $CH^{d-i}(Z)$ has $\ell$-torsion.  Then for $\ell\notin S$, the top pairing of the diagram
 is perfect.
 By Poincar{\'e} duality, the second pairing is
perfect for all $\ell\notin S$.   The commutativity of the diagram
above implies that the diagram:

$$\matrix{CH^i(Z)\otimes\bZl&\to&\mbox{Hom}(CH^{d-i}(Z)\otimes\bZl,\bZl)\cr
\downarrow&&\uparrow\cr H^{2i}(Z,\bZl(i))&\to&
\mbox{Hom}(H^{2d-2i}(Z,\bZl(d-i)),\bZl)}$$

is commutative.   Hence for $\ell\notin S$, the horizontal maps are
isomorphisms, and so the right vertical map is surjective (note that
this does not immediately follow from the injectivity of the cycle
map since the functor $\mbox{Hom}(-,\bZl)$ is not exact).
 But the right vertical map is also injective since the cokernel of the cycle map is torsion.
 Thus the left vertical map must be surjective and all maps in the
 diagram are isomorphisms.
This completes the proof of the lemma.

\begin{example}

\begin{itemize}
\item[i)] Assume that $\Yb$ is projective and that $\Yb_I$ are smooth
projective toric varieties.  Then $\Yb$ is totally degenerate.
Conditions a), b) and c) follow from ([FMSS], Corollary to Theorem
2) and ([Ku1], Proof of Theorem 6.13); see also ([D], Theorems 10.8
and 12.11), or ([F], Section 5.2, Theorem on p. 102 and the argument
on p. 103). Note that the etal{\'e} cohomology
$H^i(\Zb,\bZ/\ell\bZ)$ for $i$ odd of a toric variety $Z$ is trivial
(for example, by using the same argument as in the proof of Theorem
2.1 in [ES]), so the integral cohomology $H^i(\Zb,\bZl)$ with $i$
even is torsion free (using Theorem 2 of [FMSS]). To show that a
smooth toric variety is ordinary, one can use the fact that any
smooth proper variety $Z$ over $F$ which admits a lifting to a
smooth proper scheme over $W_2(F)$ together with the Frobenius is
ordinary. But toric varieties admit such a lifting even to $W(F)$,
as is easily seen.

\item[ii)] Let $A$ be an abelian variety over $K$. Let ${\cal A}$ be the
 N{\'e}ron model of $A$ over the ring of integers of $K$. Assume that
 the connected component of identity of the
  special fibre of ${\cal A}$ over $F$ is a split torus. We will say that
$A$ has {\it
 completely split toric reduction}. By a theorem of
 K{\"u}nnemann ([Ku1], Theorem 4.6(iii)), after passing to a finite extension
 of $K$, if necessary, $A$ has a regular projective
 semi-stable model $\tilde{\cal A}$ whose special fibre is a reduced
 divisor with smooth components, which implies that it has strictly normal
 crossings. In this case, the components of the special fibre are
 smooth projective toric varieties. The intersection of any number of
 components is also a smooth projective toric variety.  Thus $A$
 has totally degenerate reduction.

\item[iii)]Let $\widehat{\Omega}^d_K$ be the Drinfeld upper half
space obtained by removing the $K$-rational hyperplanes from ${\bf
P}^d_K$.  This is a rigid analytic space, and its quotient by a
cocompact torsion free discrete subgroup $\Gamma\subset
PGL_{d+1}(K)$ has the structure of a smooth projective variety (see
[It] for a summary of the results in this direction, and [Mus],
Theorem 4.1 and [Mu1], Theorem 2.2.5 for the original proofs).  Ito
proves the Hodge index theorem for the groups $CH^*(Y_I)_{\bQ}$ (see
[It], Conjecture 2.6 and Proposition 3.4). He also proves some of
the results in sections \S\S 3-5 of this paper. Note that in this
example, the components of the special fibre and their intersections
are successive blow-ups of projective spaces along closed linear
subschemes, and these are ordinary (see [I2], proposition 1.6).
Therefore, the results of this paper apply to these Drinfeld modular
varieties.

\item[iv)] There are Shimura varieties that can be $p$-adically uniformized by
analytic spaces other than Drinfeld upper half planes (see e.g.
[Rap], [RZ2] and [Va]). As pointed out in [La] and the introduction
to [RZ2], the bad reduction of a Shimura variety is only rarely
totally degenerate. In any case, when they can be uniformized,
Shimura varieties should have totally degenerate reduction.
Unfortunately, it seems difficult to construct the good models we
need for this paper (but see the next remark).
\end{itemize}
\end{example}

\begin{remark}  Using the results of de Jong [dJ], it should be
possible to prove that for any variety $X$ which ``looks like'' it
should have totally degenerate reduction, there is an alteration
$Y\to X$ such that $Y$ has strict semi-stable reduction and a motive
for numerical equivalence on $Y$ whose cohomology is isomorphic to
that of $X$. For example, in the theory of uniformization of Shimura
varieties (see e.g. [RZ2], [Va]), one constructs models of Shimura
varieties over $p$-adic integer rings which are often not even
regular, and it is not known at this time how to provide models such
as that required in this paper.  But the special fibres of such
models are often given by possibly singular toric varieties, and
hence the generic fibre
 should have totally degenerate reduction.
\end{remark}

\section{Definition of the Chow complexes $C^i_j(\Yb)$}

In this section we define \textit{Chow complexes} $C^i_j(\Yb)$ for
each $j=0,...,d:=\dim(Y)$. These will give us the $\bZ$-structures
we require later.  The definition is inspired by the work of Deligne
[De1], Steenbrink [St] and Guill\'en/Navarro-Aznar [GN] on mixed
Hodge theory, Rapoport and Zink on the monodromy filtration [RZ1],
Bloch, Gillet and Soul{\'e} [BGS] and Consani [Con] on the monodromy
filtration and Euler factors of L-functions.  Such
$\bZ$-structures have also been studied in some cases by Andr{\'e} [A].\\

Consider $R^{sh}$ the strict henselization of $R$, which is a
possibly ramified extension of $W$ and an unramified extension of
$R$. The assumption that $X$ has a regular proper model $\cX$ over
$R$ which is strictly semi-stable gives us  a regular proper scheme
$\cXb:=\cX \otimes_{R}R^{sh}$ over $R^{sh}$ with strictly
semi-stable reduction.\\

As in \S 1, we write $\Yb=\bigcup_{i=1}^n\Yb_i$ with each $\Yb_i$
irreducible, and for a subset $I=\{i_1,...,i_m\}$ of $\{1,...,n\}$
with $i_1<...<i_m$, denote
$$\Yb_I=\Yb_{i_1}\cap ...\cap \Yb_{i_m}.$$
By assumption, $\Yb_I$ is either empty or a reduced closed
subscheme of codimension $m$ in $\cXb$.\\

We then write $\Yb^{(m)}:=\bigsqcup \Yb_I$, where the disjoint
union is taken over all subsets $I$ of $\{1,...,n\}$ with $\# I =m$.\\

For each triple of integers $(i,j,k)$ we define
$$C^{i,k}_j:=C^{i,k}_j(\Yb)=CH^{i+j-k}(\Yb^{(2k-i+1)})$$
if $k\geq\max\{0,i\}$ and $0$ if $k<\max\{0,i\}$. Now, for each pair
of integers $(i,j)$ we define
$$C^i_j=C^i_j(\Yb)=\bigoplus_{k} C^{i,k}_j =\bigoplus_{k\geq\max\{0,i\}}
CH^{i+j-k}(\Yb^{(2k-i+1)}).$$ Note that there are only a finite
number of summands here, because $k$ runs from $\max\{i,0\}$ to
$i+j$. Note also that $C^i_j$ can be non-zero only if $i=-d,...,d$
and $j=-i,...,d-i$.

 \begin{example}
 $C^i_0=CH^0(\Yb^{(i+1)})$, $C^i_d=CH^{d+i}(\Yb^{(1-i)})$,
$C^i_{-i}=CH^0(\Yb^{(1-i)})$
 and $C^i_{d-i}=CH^{d-i}(\Yb^{(i+1)})$.
 \end{example}

\begin{observation}
For the convenience of the reader, we compare our notation with that
used in [GN] and [BGS]. Our $CH^q(\Yb_I)$ is denoted by $A^q(\Yb_I)$
in [BGS].  For the groups $K^{ij}$ used in [GN], we have

$$K^{ij}=C^i_{\frac{-i+j+d}2}$$
and
$$C^i_{j}=K^{i,i+2j-d}.$$
Note that $K^{i,j}=0$ if $-i+j+d$ is odd by definition.
\end{observation}

 Denote by
 $\rho_r:\,\Yb_{i_1...i_{m+1}} \to \Yb_{i_1...\widehat{i_r}...i_{m+1}}$
 the inclusion maps. Here, as usual, $\Yb_{i_1...\widehat{i_r}...i_{m+1}}$
 denotes the intersection of the subvarieties $\Yb_{i_j}$ for
 $j=1,...,r-1,r+1,...m+1$ (delete $\Yb_{i_r}$).\\

Now, for all $i$ and $m=1, \dots, d+1$ we have maps
$$\theta_{i,m}:\,CH^i(\Yb^{(m)})\to CH^i(\Yb^{(m+1)})$$
defined by $\theta_{i,m}=\sum_{r=1}^{m+1}(-1)^{r-1}\rho_{r}^*$,
where $\rho_r^*$ is the restriction map. Let $\rho_{r*}$ be the
Gysin map
$$\rho_{r*}:CH^i(\Yb_{i_1...i_{m+1}})\to
CH^{i+1}(\Yb_{i_1...\widehat{i_r}...i_{m+1}}),$$ and let
$$\delta_{i,m}:\,CH^i(\Yb^{(m+1)})\to CH^{i+1}(\Yb^{(m)})$$
be the map defined by $\delta_{i,m}=\sum_{r=1}^{m+1}(-1)^r\rho_{r
*}$. Let $$d'=\bigoplus_{k\geq\max\{0,i\}}\theta_{i+j-k,2k-i+1}$$
and
$$d''=\bigoplus_{k\geq\max\{0,i\}}\delta_{i+j-k,2k-i}.$$

Then define maps $d^i_j:\, C^i_j\to C^{i+1}_j$ by $d^i_j=d'+d''$.

 \begin{lemma} We have that $d^i_jd^{i-1}_j=0$, and so we get a complex for
 each $j=0,...,\mbox{dim}(\Yb)$.
 \end{lemma}

 {\bf Proof} The statement of the lemma is easily deduced from the following
 facts: For every $i$ and $m$, $\delta _{i+1,m-1} \delta _{i,m}=0$,
 $\theta _{i,m+1}\theta _{i,m}=0$ and $\theta _{i+1,m}\delta _{i,m}+ \delta
 _{i,m}\theta _{i,m+1}=0$ . The first two equalities are easy and the last
 one is proved in [BGS], Lemma 1.  Note that it is crucial for this proof that $Y$ is a reduced principal divisor with normal crossings in the
 regular scheme $\cX$.\\

 Define $T^i_j=\mbox{Ker }d^i_j/\mbox{Im }d^{i-1}_j$, the homology in
 degree $i$ of the complex $C^i_j$ defined above.\\

 We define pairings:
 $$C^{i,k}_j\times C^{-i,k-i}_{d-j}\to
C^{i,k}_{d-k}=CH^{d-(2k-i)}(\Yb^{(2k-i+1)})\stackrel{\mbox{\tiny{$(-1)^{i+j}$deg}}}{\longrightarrow}\bZ$$
 as the intersection pairing. This allows us to define
 pairings
 $$(\ ,\ ):C^i_j\times C^{-i}_{d-j}\to \bZ,$$
 by pairing each summand in $C^i_j$ with the appropriate
 summand in $C^{-i}_{d-j}$. Summands not complementary pair
 to zero. By the projection formula, these pairings are compatible with the
 differentials: $(d'x,y)=(x,d''y)$ and $(d''x,y)=(x,d'y)$. Hence they induce pairings
  $$(\ ,\ ): T^i_j\times T^{-i}_{d-j}\to \bZ.$$

The monodromy operator $N:\,C^i_j\to C^{i+2}_{j-1}$ is defined as
the identity map on the summands in common, and the zero map on
different summands (observe that $C^{i,k}_j$ and $C^{i+2,k+1}_{j-1}$
are equal). Clearly, $N$ commutes with the differentials, and so
induces an operator on the $T^i_j$, which we also denote by $N$. We
have that $N^i$ is the identity on $C^i_j$ for $i\geq 0$. This is
proved by showing that the summands in a given $C^i_j$ persist
throughout, and the others occurring in subsequent groups
eventually disappear.\\

The following result is a direct consequence of a result of Guillen
and Navarro ([GN] Prop. 2.9 and Th{\'e}or{\`e}me 5.2 or [BGS], Lemma
1.5 and Theorem 2), using the fact that the Chow groups of the
components of $\Yb$ in the case of totally degenerate reduction
satisfy the hard Lefschetz theorem and the Hodge index theorem. It
is a crucial result for this paper.

 \begin{proposition} Suppose that $\Yb$ satisfies the conditions
 in Definition 1. Then the ith power $N^i$ of the monodromy operator $N$ induces an isogeny:
 $$N^i:\,T^{-i}_{j+i}\to T^{i}_{j}$$
 for all $i\geq 0$ and $j$. Moreover the pairings
  $$(\ ,\ ): T^i_j\times T^{-i}_{d-j}\to \bZ$$
are nondegenerate on the torsion free quotients.
 \end{proposition}

\section{The monodromy filtration in the case $\ell\ne p$}

In this section, we study the monodromy filtration on the {\'e}tale
cohomology $H^*(\Xb,\bQl)$ using the techniques developed in  \S 3.
The goal is to establish an isomorphism between the graded quotients
for the monodromy filtration and an appropriate twist of
$T^i_j\otimes_{\bZ}\bQl$ for any prime number $l\ne p$.
In \S 5 and \S 6, we will consider the case $\ell=p$.\\

Recall Grothendieck's monodromy theorem in this situation (see
[SGA\,VII, I] or [ST], Appendix for the statement, and [RZ1] for the
proof): the restriction to the inertia group $I$ of the $\ell$-adic
representation associated to the {\'e}tale cohomology
$H^i(\Xb,\bQl)$ of a smooth projective variety $X$ with semistable
reduction over a complete discretely
valued field $K$ is unipotent.\\

Thus we have the monodromy operator $N:H^i(\Xb,\bQl)(1) \to
H^i(\Xb,\bQl)$, characterized by the fact that the restriction of
the representation to the maximal pro-$\ell$ quotient of the inertia
group $I$ is defined by the composition $\exp \circ N\circ
t_{\ell}$, where $t_{\ell}:I\to \bZl(1)$ is the natural map. Using
this monodromy operator, we can construct ([De3], 1.6.1, 1.6.13) the
monodromy filtration, $M_{\bullet}$, which is the unique increasing
filtration
$$0=M_{-i-1}\subseteq M_{-i}\subseteq...\subseteq M_{i-1}\subseteq
M_i=H^i(\Xb,\bQl),$$ such that
$$ NM_j(1)\subseteq M_{j-2}$$
and
$$ N^j\, \mbox{induces an isomorphism}\,
\mbox{Gr}_j^M(j)\cong \mbox{Gr}_{-j}^M.$$

The construction of the monodromy filtration in ibid. shows that $I$
acts trivially on the graded quotients $\mbox{Gr}_j^M$ and
$M^I\subseteq M_0$.  If the weights of Frobenius are integers, then
we also have the weight filtration ([De3],
Proposition-D{\'e}finition 1.7.5): the unique increasing filtration
on $H^i(\Xb,\bQl)$ whose graded quotients are pure
$Gal(\bF/F)$-modules. Recall also the mo\-no\-dro\-my-weight
conjecture ([De1], Principe 8.1; [RZ], Introduction), which can be
stated in two equivalent forms:

 \begin{enumerate}
 \item
 $\mbox{Gr}_j^M$ is a pure
$\mbox{Gal}(\bF/F)$-module of weight $i+j$.

\item The weight
filtration has the properties mentioned above that characterize the
monodromy filtration. \end{enumerate}

This conjecture in its first form was proved by Deligne in the
equi-cha\-racte\-ris\-tic case ([De2], Th{\'e}or{\`e}me 1.8.4), and
in its second form by Rapoport and Zink in the case where $X$ is of
dimension 2 and has a model $\cX$ with special fibre $Y=\sum n_iY_i$
(as cycle), with $n_i$ prime to $l$ for each $i$ (see
[RZ], Satz 2.13).   We will use the second form below.\\

Observe that when $\ell=p$, there is in general no ``simple''
monodromy filtration, as is pointed out in ([Ja1], page 345). That
is the reason for treating this case separately in the next two
sections.\\

Consider the weight spectral sequence of Rapoport-Zink ([RZ], Satz
2.10; see also [Ja2]),  in its two forms:
$$E_1^{i,j}=H^{i+j}(\Yb,\mbox{Gr}^W_{-i}(R\Psi(\bZl))\Longrightarrow
H^{i+j}(\Xb,\bZl).$$ and
$$E_1^{i,j}=\bigoplus_{k\geq\max\{0,i\}}
H^{j+2i-2k}(\Yb^{(2k+1-i)},\bZl(i-k))\Longrightarrow
H^{i+j}(\Xb,\bZl).$$ By Deligne's weight purity theorem (Riemann
hypothesis for smooth projective varieties over finite fields), this
spectral sequence tensored by $\bQl$ degenerates at $E_2$ when the
residue field is finite, and Nakayama proved the degeneration in the
general case [Na]. We denote by $\free{E^{i,j}_{r}}$ the
quotient of $E_r^{i,j}$ by its torsion subgroup.\\

        We write the two forms of this spectral sequence because the
first one has simple indexing, and the groups in the second are very
similar to the ones considered in the last section, as we now see:

\begin{proposition}  Suppose that $Y$ satisfies the condition b) in Definition 1.
Then we have:
$$\free{E^{i,2j}_{1}}\cong \free{C^i_j(\Yb)}\otimes \bZl(-j)$$
$$\free{E^{i,2j+1}_{1}}=0.$$
These are compatible with the differentials on both sides, the
differentials on the left hand side being those defined in section
3, tensored with $\bZl$.
\end{proposition}

{\bf Proof}: The two equalities are clear from our assumptions on
$\Yb$. The compatibility between differentials is deduced directly
from the results of Rapoport and Zink ([RZ], \S 2, especially Satz
2.10; see also [Ja2], \S 3).\\

Let us denote by $W_{\bullet}$ the filtration on $H^{n}(\Xb,\bZl)$
induced by the weight spectral sequence, and by $M_{\bullet}$ the
filtration defined by $M_i(H^n):=W_{i+n}(H^n)$. Then the filtration
$W_{\bullet}$ gives us the weight filtration on $H^{n}(\Xb,\bQl)$,
so the monodromy conjecture says that the filtration $M_{\bullet}$
is the monodromy filtration.

\begin{corollary}  With assumptions as in Proposition 2, we have canonical
isogenies compatible with the action of the Galois group $G$
$$\free{T^i_j(\Yb)}\otimes\bZl(-j)\to
\free{\mbox{Gr}^M_{-i}H^{i+2j}(\Xb,\bZl)},$$ and
$\mbox{Gr}^M_{-i}H^{i+2j+1}(\Xb,\bZl)$ is torsion.

Moreover, if $Y$ verifies condition a) in Definition 1, these
isogenies are isomorphisms for almost all $\ell$.
\end{corollary}

{\bf Proof}: By  Theorem 0.1 in [Na], we know that the weight
spectral sequence degenerates at $E_2$ after tensoring by $\bQl$.
This implies that the natural maps
$$\mbox{Gr}^M_{-i}H^{i+2j}(\Xb,\bZl)=\mbox{Gr}^W_{2j}H^{i+2j}(\Xb,\bZl)\to
E^{i,2j}_{2}\cong T^i_j(\Yb)$$ are isomorphisms after tensoring by
$\bQl$ and so the induced maps on the torsion free quotients are
isogenies. It also implies that
$\mbox{Gr}^M_{-i}H^{i+2j+1}(\Xb,\bQl)=0$.

To see the second assertion, recall that under the assumption the
groups $T^i_j(\Yb)$ are finitely generated abelian groups, so
$T^i_j(\Yb)\otimes\bZl$ are all torsion free for almost all $\ell$.
For these primes $\ell$, we can deduce the degeneration at $E_2$
from the degeneration after tensoring by $\bQl$.

\begin{corollary} Supose that $Y$ satisfies conditions a) and b) in Definition 1.
Then the filtration $M_{\bullet}\otimes_{\bZl}\bQl$ is the monodromy
filtration and the induced maps
$$\mbox{Gr}^M_{-i}H^{i+2j}(\Xb,\bQl)
\stackrel{N}\rightarrow \mbox{Gr}^M_{-i-2}H^{i+2j}(\Xb,\bQl)(1)$$
coincide under the isomorphism in Corollary 1 with the monodromy
operator defined in \S 3.
\end{corollary}
{\bf Proof}: It is not difficult to see that $NM_i(1)\subseteq
M_{i-2}$. Thus one only has to show that the $i$-th power of the
monodromy induces an isomorphism between $\mbox{Gr}_{i}^M H^n \cong
\mbox{Gr}_{-i}^M H^n(-i).$ Using the corollary above, this is
reduced to showing that $N^i$ induces an isomorphism
$$\,T^{-i}_{j+i}\otimes \bQl(-j-i) \to T^{i}_{j}\otimes \bQl(-j-i).$$
In ([RZ], \S 2), it is shown that the map
$$T^{i}_{j}\otimes\bQl(-j)\!\cong\!\mbox{Gr}^M_{-i}H^{i+2j}(\Xb,\bQl)
\!\stackrel{N}\rightarrow \!\mbox{Gr}^M_{-i-2}H^{i+2j}(\Xb,\bQl)(-1)
\!\cong \! T^{i+2}_{j-1}\otimes\bQl(-j)$$ is the same as the map
that we constructed in section 2, and hence we can deduce the result
from Proposition 1.

\begin{remark}
Observe that for almost all $\ell$, we have an isomorphism
$$T^i_j(\Yb)\otimes\bZl(-j)\cong
\mbox{Gr}^M_{-i}H^{i+2j}(\Xb,\bZl)$$ and hence the graded quotients
of the monodromy filtration are torsion free for almost all $\ell$.

This is because by lemma 1 we have
$$E^{i,2j}_{1}\cong C^i_j(\Yb)\otimes \bZl(-j)$$
$$E^{i,2j+1}_{1}=0$$
for almost all $\ell$, and using the compatibility of the
differentials, we get that
$$E^{i,2j}_{2}\cong T^i_j(\Yb)\otimes \bZl(-j)$$
$$E^{i,2j+1}_{2}=0.$$
so the $E_2$-terms are torsion free for almost all $\ell$ because
the groups $T^i_j$ are finitely generated abelian groups. For such
$\ell$, the proof of Corollary 1 shows that the Rapoport-Zink
spectral sequence degenerates at $E_2$.
\end{remark}

\begin{remark}  In papers of Ito [It] and de Shalit [dS], the monodromy and weight
filtrations are considered in $\ell$-adic and log-crystalline
cohomology.  The appropriate versions of the monodromy-weight
conjecture are proved for varieties uniformized by the Drinfeld
upper half space. We compare and contrast our results with theirs.
Whereas we assume the Hodge index theorem for the Chow groups of
intersections of components of the special fibre of a regular proper
model of our $X$, Ito proves this for the varieties he considers.
Thus Ito's paper provides us with more examples to which the methods
of this paper may be applied. de Shalit proves the $p$-adic version
of the monodromy-weight conjecture by doing harmonic analysis and
combinatorics on the Bruhat-Tits building of $PGL_{d+1}(K)$.  On the
other hand, their results and methods do not give a monodromy
filtration on the $p$-adic cohomology of $\Xb$, as we do here.
\end{remark}

\section{The monodromy filtration on log-crystalline cohomology}

Recall the notations we introduced in \S 1: $F$ is a perfect field
of characteristic $p>0$,  $W(F)$ is the ring of Witt vectors with
coefficients in $F$, and $K_0$ is the field of fractions of $W(F)$.
We will denote $W(\Fb)$ by $W$ and denote its fraction
field by $L$.\\

Let ${\cal L}$ be a logarithmic structure on $F$, in the sense of
Fontaine, Illusie and Kato (see [Ka]). Let $(Y,{\cal M})$ be a
smooth $(F,{\cal L})$-log-scheme, i.e. $Y$ is a scheme over $F$ and
${\cal M}$ is a fine log-structure with a smooth map of log schemes
$f:(Y,{\cal M}) \to (F,{\cal L})$. We assume also that $(Y,{\cal
M})$ is semistable in the sense of ([Mok], D{\'e}finition 2.4.1).
Hyodo and Kato [HK] have defined under these conditions the
log-crystalline cohomology of $Y$, which we will denote by
$H^{i}(Y^{\times}/W(F)^{\times})$. It is a $W(F)$-module of finite
type, with a Frobenius operator $\Phi $ that is bijective after
tensoring by $K_0$ and semilinear with respect to the Frobenius in
$W(F)$, and a monodromy operator $N$ that is nilpotent and satisfies
$N\Phi=p\Phi N$. This operator $N$ then  determines a filtration on
$H^{i}(Y^{\times}/W(F)^{\times})$ that we call the {\it monodromy
filtration}. Recall that the log-crystalline cohomology can be
computed as the inverse limit of the hypercohomology of the de
Rham-Witt complex
$W_n\omega^{\bullet} _Y$ of level $n$.\\

Suppose now that $Y$ is a proper variety. Then we have the spectral
sequence of Mokrane ([Mok], 3.23):
$$E_1^{i,j}(Y):=\! \! \! \! \bigoplus_{k\geq\max\{0,i\}} \! \! \! \!
 H^{j+2i-2k}(Y^{(2k+1-i)}/W(F))(i-k) \Longrightarrow
 H^{i+j}(Y^{\times}/W(F)^{\times}),$$
where the twist by $i-k$ means to multiply the Frobenius operator by
$p^{k-i}$. We denote by $\free{E^{i,j}_{1}}$ the free quotient of
the terms of this spectral sequence.  The associated filtration on
the abutment is called the {\it weight filtration}, and as before we
will denote it by $W_{\bullet}$. We will also denote by
$M_{\bullet}$ the filtration defined by $M_i(H^n):=W_{i+n}(H^n)$.

We can consider also the log-crystalline cohomology
$H^{n}(\Yb^{\times}/W^{\times})$ of $\Yb$, which is a $W$-module of
finite type, with a Frobenius operator $\Phi $ and a monodromy
operator, and which  also has an action of the absolute Galois group
$Gal(\Fb/F)$. We then have a canonical isomorphism as $W(F)$-modules
with a Frobenius
$$H^{n}(Y^{\times}/W(F)^{\times}) \cong
H^{n}(\Yb^{\times}/W^{\times})^{Gal(\Fb/F)}.$$ Moreover, Mokrane's
spectral sequence
$$E_1^{i,j}(\Yb):=\bigoplus_{k\geq\max\{0,i\}}
 H^{j+2i-2k}(\Yb^{(2k+1-i)}/W)(i-k) \Longrightarrow
 H^{i+j}(\Yb^{\times}/W^{\times})$$
is compatible also with the Galois action and the spectral sequence
for $Y$ above is canonically isomorphic to the
$Gal(\Fb/F)$-invariant part of this spectral sequence.

  \begin{proposition} With notation and assumptions as above,
  suppose further that $Y$ satisfies the assumption c) of Definition 1.
  Then we have a canonical isomorphism compatible with the action of
  $Gal(\Fb/F)$ and with the Frobenius
  $$\free{E^{i,2j}_{1}(\Yb)}\cong \free{C^i_j(\Yb)}\otimes W(-j)$$
  $$\free{E^{i,2j+1}_{1}(\Yb)}=0.$$
  These are compatible with the differentials on both sides, the
  differentials on the right hand side being those defined in \S 3.

  We also have that $E^{i,2j}_{1}(Y)(j)\otimes_{W(F)}K_0$ corresponds to the $p$-adic
  representation $C^i_j(\Yb)\otimes \bZp$ of $Gal(\Fb/F)$ by the
  correspondence explained in Remark 1 (ii), and
  $E^{i,2j+1}_{1}(Y)\otimes_{W(F)}K_0=0$.
  \end{proposition}

{\bf Proof}: The first two equalities are clear by condition c) of
\S 1, and the compatibility between the differentials is clear by
using the results of Mokrane ([Mok], \S 3), who proves that the
spectral sequence degenerates at $E_2$ modulo torsion for any $Y$
([Mok], 3.32(2)). For the last assertion, apply Remark 1, (ii) of \S
2.

  \begin{corollary} Suppose that $Y$ satisfies a) and c) in definition 1.
  Then the spectral sequence of Mokrane degenerates at $E_2$,
  the weight filtration induces the monodromy
  filtration, and there is a canonical isogeny compatible with the
  actions of $Gal(\Fb/F)$ and the Frobenius
  $$\free{T^i_j(\Yb)}\otimes W(-j)\to
  \free{\mbox{Gr}^M_{-i}H^{i+2j}(\Yb^{\times}/W^{\times})},$$
  and $\mbox{Gr}^M_{-i}H^{i+2j+1}(\Yb^{\times}/W^{\times})$ is
  torsion.
  Moreover, the monodromy map on the graded quotients agrees with
  the map defined in \S 3.

  We also have that the $\varphi$-module
  $\mbox{Gr}^M_{-i}H^{i+2j}(Y^{\times}/W(F)^{\times})\otimes_{W(F)}K_0(j)$ corresponds
  to the p-adic representation $T^i_j(\Yb)\otimes \bQp$ of $Gal(\Fb/F)$.
  \end{corollary}

{\bf Proof} The degeneration of the spectral sequence is deduced
using a slope argument. Observe that for any $i$ and any $k$, the
cohomology groups
  $$\free{H^{2k}(\Yb^{(i)}/W)(k)}$$
are of pure slope $0$ because they are generated by algebraic cycles
by our conditions in \S 2, c). So, the $E_1^{i,2j}$-terms are of
pure slope $j$ modulo torsion, as are any of their subquotients,
such as the $E_2^{i,2j}$ terms. So, we have that
$\mbox{Hom}(E_2^{r,s},E_2^{r+t,s-t+1})=0$ modulo torsion for all $r$
and $s$ and $t\geq 2$, because these are of different slopes. Hence
$\free{E_2^{r,s}}=\free{E_t^{r,s}}$ for all $r,s$ and $t\geq 2$, and
by an easy induction the differentials are zero on the
$E_t^{r,s}$-terms modulo torsion for all $r,s$ and $t\geq 2$.

The assertion about the coincidence of the weight filtration and the
monodromy filtration is deduced from Proposition 1 by using ([Mok],
Proposition 3.18) as in the proof of Corollary 2.

The other assertions are deduced easily from Proposition 3,
using Remark 1, (ii) of \S 2.\\

Now, also assume that $Y$ is ordinary (see condition d) of \S 1).
Note that this condition is stable under base change to $\Fb$. Then
there exists a canonical decomposition as $W$-modules with a
Frobenius (see [I3], Proposition 1.5. (b))
 $$ H^{n}(Y^{\times}/W^{\times}) \cong \bigoplus_{i+j=n} H^j(Y,W\omega
^i)(-i).$$
We will identify the two sides of the isomorphism.\\

The existence of this canonical decomposition allows us to define
the filtration by \textit{increasing slopes} $U_j$ as
 $$U_j(H^{n}(Y^{\times}/W^{\times})):=\oplus_{r\leq j}
 H^{n-r}(Y,W\omega^r)(-r),$$
which has the property that
 $$ \mbox{Gr}^U_{j}H^{n}(Y^{\times}/W^{\times})\cong H^{n-j}(Y,W\omega
^j)(-j).$$

This filtration is \textit{opposed} to the usual \textit{slope
filtration}, which is given by decreasing slopes.  We can do the
same thing for $\Yb$.

  \begin{corollary} Suppose that $Y$ satisfies c) and d) in Definition 1.
  We have then that
  $$M_{2j-n}(\free{H^{n}(\Yb^{\times}/W^{\times})})\subseteq
  M_{2j-n+1}(\free{H^{n}(\Yb^{\times}/W^{\times})}) $$
  $$\lefteqn{ \subseteq U_j(\free{H^{n}(\Yb^{\times}/W^{\times})})},$$
  with cokernel of finite exponent. Thus if we tensor with $L$  we have equalities,
  and this holds also when $\Yb$ is replaced by $Y$ and we tensor by $K_0$.\\

  As a consequence, if $Y$ also satisfies condition a) in Definition 1,
  we have canonical isogenies
  $$ \free{T^{i-j}_j(\Yb)}\otimes W \to \free{H^{i}(\Yb, W\omega^j_{\Yb})}.$$
 \end{corollary}

  {\bf Proof.} Observe first that the filtration by increasing slopes has a
  splitting and hence the torsion free part of
  $\mbox{Gr}^U_{j}H^{i+2j}(\Yb^{\times}/W^{\times})$
  is equal to
  $\mbox{Gr}^U_{j}\left(\free{H^{i+2j}(\Yb^{\times}/W^{\times})}\right)$.

  By Corollary 3, we have that
  $\free{\mbox{Gr}^M_{-i}H^{i+2j}(\Yb^{\times}/W^{\times})}$ is of
  pure slope $j$, and
  $\mbox{Gr}^M_{-i}H^{i+2j+1}(\Yb^{\times}/W^{\times})$ is
  torsion, hence the first assertion is true.

  From these facts, we have that $\free{\mbox{Gr}^M_{-i}H^{i+2j}(\Yb^{\times}/W^{\times})}$
  is canonically isogenous to
  $$\free{\mbox{Gr}^U_{j}H^{i+2j}(\Yb^{\times}/W^{\times})}\cong
  \free{H^{i+j}(\Yb,W\omega^j)}(-j).$$
  Composing with the canonical isogeny
  $$\free{T^i_j(\Yb)}\otimes W(-j)\to
  \free{\mbox{Gr}^M_{-i}H^{i+2j}(\Yb^{\times}/W^{\times})},$$
  we get a canonical isogeny
  $$\free{T^i_j(\Yb)}\otimes W\to \free{H^{i+j}(\Yb,W\omega^j)}.$$\\

Now let's consider the logarithmic Hodge-Witt pro-sheaves
$W\omega^j_{\Yb ,log}$ on $\Yb$, defined as the kernel of $1-F$ on
the pro-sheaves $W\omega^j_{\Yb}$ for the {\'e}tale site.\\

To prove the next result we will use Proposition 2.3 in [I3], which
says that
$$H^i(\Yb, W\omega^j_{\Yb,log})\otimes_{\bZp} W \cong
 H^i(\Yb, W\omega^j_{\Yb}).$$

\begin{corollary} Suppose that $Y$ satisfies a), c) and d) in Definition 1.
Then we have a canonical isogeny
$$\free{T^{i-j}_j(\Yb)}\otimes_{\bZ}{\bZp} \to
\free{H^{i}(\Yb,W\omega^j_{\Yb,log})}.$$
\end{corollary}

{\bf Proof.} By Corollary 4 we have a canonical isogeny
$$ T^{i-j}_j(\Yb) \otimes_{\bZ} W(j) \to  H^i(\Yb,
W\omega^j_{\Yb}) \otimes_{\bZp} W$$ compatible with the actions of
$Gal(\Fb/F)$ and the Frobenius automorphisms. Here $(j)$ means the
twist of Frobenius (not the Galois action). By taking the fixed part
$Frob=p^j$ of both sides, we obtain a canonical isogeny
$$ T^{i-j}_j(\Yb) \otimes_{\bZ} \bZp \to H^i(\Yb,
W\omega^j_{\Yb,log})$$ compatible with the actions of $Gal(\Fb,F)$,
using the fact that the action of the Frobenius on
$H^i(\Yb,W\omega^j_{\Yb,log})$ is multiplication
by $p^j$.\\

  \section{The monodromy filtration on $p$-adic cohomology}

The goal of this section is to show that the filtration on the
$p$-adic cohomology of $X$ induced by the vanishing cycles spectral
sequence has all of the properties that a monodromy filtration
should have when $X$ has totally degenerate
reduction.\\

First of all, recall that a $p$-adic representation $V$ of $G_K$ is
{\it ordinary} if there exists a filtration $(Fil ^i V)_{i\in \bZ}$
of $V$, stable by the action of $G_K$, such that the inertia
subgroup $I_K$ acts on $Fil ^i V/Fil ^{i+1} V$ via $\chi ^i$, where
$\chi $ denotes the cyclotomic character. It is easy to see that
this filtration is unique. The next theorem is the main result in
[H] (see also [I3], Th{\'e}or{\`e}me 2.5 and Corollaire 2.7 ).

 \begin{theorem} (Hyodo) Assume that $Y$ is ordinary.
 Then the $p$-adic representation $H^m(\Xb,\bQp)$ is
 ordinary. Moreover, the vanishing cycles spectral sequence
 $$E_2^{i,j}=H^{i}(\Yb,R^{j}\Psi(\bZp))\Longrightarrow
 H^{i+j}(\Xb,\bZp)$$
 degenerates at $E_2$ modulo torsion, and if $F^{\bullet}$ denotes the
corresponding
 filtration
 on $H^{i+j}(\Xb,\bZp)$, one has a canonical isogenies as
 $G_K$-modules
 $$\free{\mbox{Gr}_F^{i}H^{i+j}(\Xb,\bZp)}\to
 \free{H^{i}(\Yb, W\omega^j_{\Yb,log})}(-j).$$
 \end{theorem}

 Suppose now that $\cX$ has special fibre that satisfies the assumptions
of \S 1 and \S 2. Then, as a consequence of Corollary 5 in \S 5, we
get the following result.

\begin{corollary}
 Assume that $Y$ satisfies a), c) and d) in Definition 1. Then, if
 $F^{\bullet}$ denotes the filtration obtained from the vanishing
 cycles spectral sequence, one has canonical isogenies as
 $G_K$-modules
 $$\free{T^{i}_j(\Yb)}\otimes \bZp(-j)\to
 \free{\mbox{Gr}_F^{i+j}H^{i+2j}(\Xb,\bZp)}$$
\end{corollary}

We have then that the filtration $F^{\bullet}$ has the same type of
graded quotients as the monodromy filtration in $\ell$-adic
cohomology (modulo isogeny). In fact we will prove that this
filtration can be obtained from the monodromy filtration on the
log-crystalline cohomology by
taking the functor $V_{st}$.\\

Recall that a filtered $(\Phi,N)$-module $H$ is a $K_0$-vector space
of finite dimension $H$ with a Frobenius $\Phi$ which is a
semi-linear automorphism and a monodromy map $N$ which is a
$K_0$-linear endomorphism verifying that $N\Phi=p\Phi N$ and an
decreasing filtration $Fil^{\bullet}$ in $H\otimes_{K_0}K$ by
$K$-subspaces which is exhaustive and separated. Given a filtered
$(\Phi,N)$-module $H$, we define
$$V_{st}(H):=(B_{st}\otimes_{K_0}
H)^{N=0, \Phi=1}\cap Fil^0(B_{dR}\otimes_{K} (H\otimes_{K_0} K)).$$
It is a $p$-adic representation of $G_K$, that is a $\bQp$-vector
space with a continuous action of $G_K$. Recall from \S 1 that we
fixed a choice of uniformizer $\pi$ of $K$, which determined the
embedding of $B_{st}$ in $B_{DR}$, and so this
$p$-adic representation depends on that choice.\\

Consider the log-crystalline cohomology $H^n(Y^{\times}/W^{\times})$
we used in the last section. Hyodo and Kato proved in [HK], \S 5,
that we have an isomorphism
$$\rho_{\pi}\colon  H^n(Y^{\times}/W^{\times})\otimes_{W}K \cong
H^n_{dR}(X/K)$$ depending on the uniformizer $\pi$ that we have
chosen. Using this isomorphism, we get a structure of filtered
$(\Phi,N)$-module on
the log-crystalline cohomology of the log-scheme $Y$.\\

\textit{Now assume only that $Y$ is ordinary}, and consider the
filtration by increasing slopes $U_{\bullet}$ that we discussed in
\S 5. Then the induced $(\Phi,N)$-module structure on the $U_i$ via
the Hyodo-Kato isomorphism $\rho_{\pi}$ gives a filtration by
filtered $(\Phi,N)$-modules. Applying the functor $V_{st}$, and
using the main result of Tsuji ([Ts], p. 235) on the conjecture
$C_{st}$
$$ H^{n}(\Xb,\bQp) \cong V_{st}(H^n(Y^{\times}/W^{\times})),$$
we get a filtration on the $p$-adic cohomology $H^n(\Xb,\bQp)$.

\begin{theorem}
Assume only that $Y$ is ordinary.  Then the filtrations
$V_{st}(U_i)$ and $F^{n-i}$ on $H^n(\Xb,\bQp)$ are the same.

\end{theorem}

  {\bf Proof.}  First of all, note that the filtered
  $(\Phi,N)$-module $H^n(Y^{\times}/W^{\times})$ is ordinary in the sense of ([P], p. 186).
  This is because $H^n(Y^{\times}/W^{\times})$ is isomorphic, by Tsuji's
  theorem, to $D_{st}( H^{n}(\Xb,\bQp))$, and the functor $D_{st}$ takes
  ordinary $p$-adic representations to ordinary filtered $(\Phi,N)$-modules
  (see [P], Th{\'e}or{\`e}me 1.5).\\

  Now observe that the graded quotients $D$ with respect to the $U_{\bullet}$
  filtration have the property that there exists $i$ such that
  $Fil^{i}(D_K)=D_K$
  and $Fil^{i+1}(D_K)=0$, $p^{-i}\Phi$ acts as an automorphism on a lattice in
  $D$ and $N=0$. This is due
  to the fact that the $U_{\bullet}$ filtration is opposed to the
  Hodge filtration, because our filtered $(\Phi,N)$-module is ordinary
  (see [P], middle of p. 187 and [I3] 2.6 c), middle of p. 217). But by
  ([P], Lemme 2.3), the inertia group acts on $V_{st}(D)$ as $\chi ^i$.
  Moreover, the functor $V_{st}$ is exact for the ordinary
  $(\Phi,N)$-filtered modules (see [P], 2.7).
  So the $V_{st}(U_{\bullet})$ filtration has the same graded quotients as
  the filtration $F^{n-\bullet}$, and hence these filtrations must be the
  same.\\

\begin{remark}
\begin{itemize}
\item[(i)] For $X$ with ordinary reduction, there is another
structure of $(\Phi,N)$-filtered module on the de Rham cohomology
$H^n_{dR}(X/K)$ given by the result of Hyodo [H] (see also [I3],
Corollaire 2.6 (c)). This structure is not the right one for
applying Tsuji's comparison theorem, although the filtration by
increasing slopes we get (or, equivalently, the monodromy
filtration) is also opposed to the Hodge filtration and so they have
all isomorphic graded quotients. \item[(ii)] We have canonical
isomorphisms
$$T^{i-j}_j(\Yb)\otimes_{\bZ}K\cong H^i(X,\Omega^j_X)$$
because the filtration by increasing slopes is opposed to the Hodge
filtration. Now, proposition 1 implies that the monodromy map
$N\colon T^i_j \to T^{i+2}_{j-1}$ is, after tensoring by $\bQ$,
injective if $i<0$ and surjective if $i\ge 0$. This implies that we
have a monodromy map
$$N \colon H^i(X,\Omega^j_X) \to H^{i+1}(X,\Omega^{j-1}_X)$$
which is injective if $i<j$ and surjective if $i\ge j$. So, a
necessary condition for a variety to have totally degenerate
reduction is that the dimensions $h^{j,i}:=\dim_K H^i(X,\Omega^j_X)$
satisfy that, if $n=2i$,  $$h^{n,0}\le h^{n-1,1}\le \cdots \le
h^{i,i}$$ and if $n=2i-1$,
$$h^{n,0}\le h^{n-1,1}\le \cdots \le h^{i,i-1}$$
(recall that $h^{j,i}=h^{i,j}$ by Hodge theory). This condition
excludes rigid Calabi-Yau threefolds, for example, which have
$h^{3,0}=1$ and $h^{2,1}=0$.
\end{itemize}
\end{remark}

We summarize our discussion in the following:

\begin{theorem}  For $X$ with totally degenerate reduction,
 consider the filtration $M_{\bullet}$  on $H^n=H^n(\Xb,\bQp)$ defined
 as $M_{i}(H^n):=F^{j}(H^n)$, where $j=\frac {n-i}2$ or $j=\frac {n-i+1}2$
depending
 on the parity of $n-i$. Then we have
 $$T^i_j(\Yb)\otimes\bQp(-j)\cong
 \mbox{Gr}^M_{-i}H^{i+2j}(\Xb,\bQp),$$
 and the maps on these graded quotients induced by $N$ on the
 $T^i_j$ coincide after applying the functor $V_{st}$ with the
 maps $N$ defined in the section \S 3. Moreover, they
 verify that $N^i$ induces an isomorphism $\mbox{Gr}_{i}^M \cong
 \mbox{Gr}_{-i}^M(-i)$. The other graded quotients are zero.
\end{theorem}

\centerline{\bf References}

\begin{itemize}
\item[{\bf [A]}]  Y. Andr{\'e}, $p$-adic Betti lattices, in $p$-adic
analysis, Trento 1989, 23-63, Lecture Notes in Mathematics 1454,
1990
\item[{\bf [BK]}] S. Bloch and K. Kato, $p$-adic {\'e}tale
cohomology, Publ. Math. I.H.E.S. 63 (1986) 107-152
\item[{\bf [BGS]}] S. Bloch, H. Gillet and C. Soul{\'e}, Algebraic cycles
on degenerate fibers, in Arithmetic Aspects of Algebraic Geometry,
Cortona 1994, F. Ca\-ta\-ne\-se editor, 45-69.
\item[{\bf [Con]}]  C. Consani, Double complexes and Euler
factors, Comp. Math. 111 (1998) 323-358
\item[{\bf [D]}]  V.I. Danilov, The geometry of toric varieties (in
Russian), Uspekhi Mat. Nauk 33:2 (1978), 85-134; English translation
in Russian Mathematical Surveys 33:2 (1978) 97-154.
\item[{\bf [DeJ]}]  A.J. de Jong, Smoothness, alterations and
semi-stability, Publ. Math. I.H.E.S. 83 (1996) 52-96.
\item[{\bf [De1]}]  P. Deligne, Th{\'e}orie de Hodge I, Actes du
Congr{\`e}s International des Math{\'e}maticiens, Nice, 1970, Tome
I, 425-430.
\item[{\bf [De2]}]  P. Deligne, La conjecture de Weil I, Publ. Math.
I.H.E.S. 43 (1974) 273-308.
\item[{\bf [De3]}]  P. Deligne, La conjecture de Weil II, Publ. Math.
I.H.E.S. 52 (1980) 137-252.
\item[{\bf [dS]}]  E. de Shalit, The $p$-adic monodromy-weight
conjecture for $p$-adically uniformized varieties,  Compos. Math.
141  (2005),  no. 1, 101--120.
\item[{\bf [ES]}]  G. Ellingsrud, S. A. Str{\o}mme, Towards
the Chow ring of the Hilbert scheme of ${\bf P}^2$, J. Reine Angew.
Math. 441  (1993), 33-44.
\item[{\bf [F]}]  W. Fulton, Introduction to Toric Varieties, Annals of
Math. Studies, Volume 131, Princeton University Press, 1993.
\item[{\bf [FMSS]}]  W. Fulton, R. MacPherson, F. Sottile and B.
Sturmfels, Intersection theory on spherical varieties, J. Alg.
Geometry 4 (1995) 181-193
\item[{\bf [G]}]  O. Gabber, Sur la torsion dans la cohomologie
$\ell$-adique d'une vari{\'e}t{\'e}, C.R.A.S. 297 (1983) 179-182.
\item[{\bf [GN]}] F. Guill{\'e}n and V. Navarro-Aznar, Sur le th{\'e}or{\`e}me
local des cycles invariants, Duke Math. J. 61 (1990), 133-155.
\item[{\bf [H]}] O. Hyodo, A note on $p$-adic {\'e}tale cohomology in the
semi-stable reduction case, Inv. Math. 91 (1988), 543-557.
\item[{\bf [HK]}] O. Hyodo and K. Kato, Semi-stable reduction and
crystalline cohomology with logarithmic poles, in {\it P{\'e}riodes
p-adiques, S{\'e}minarie de Bures, 1988} (J.-M. Fontaine ed.),
Ast{\'e}risque 223 (1994), 221-268.
\item[{\bf [I1]}] L. Illusie, Cohomologie de de Rham et cohomologie
{\'e}tale $p$-adique (D'apr{\`e}s G. Faltings, J.-M. Fontaine et
al.) in {\it S{\'e}minaire Bourbaki}, 1989-1990, Asterisque 189-190
(1990), 325-374
\item[{\bf [I2]}]  L. Illusie, Ordinarit{\'e} des intersections
compl{\`e}tes g{\'e}n{\'e}rales, in {\it The Grothendieck
Festschrift, Vol II}, Progr. Math, Vol. 87, Birkh{\"a}user, Boston,
1990
\item[{\bf [I3]}] L. Illusie, R{\'e}duction semi-stable ordinaire,
cohomologie {\'e}tale $p$-adique et cohomologie de de Rham,
d'apr{\`e}s Bloch-Kato et Hyodo, in {\it P{\'e}riodes  p-adiques,
S{\'e}minaire de Bures, 1988} (J.-M. Fontaine ed.), Ast{\'e}risque
223 (1994), 209-220.
\item[{\bf [I4]}]  L. Illusie, Autour du th{\'e}or{\`e}me de
monodromie locale, in {\it P{\'e}riodes $p$-adiques, S{\'e}minaire
de Bures, 1988} (J.-M. Fontaine ed.), Ast{\'e}risque 223 (1994),
9-57.
\item[{\bf [It]}]  T. Ito, Weight-Monodromy conjecture for
$p$-adically uniformized varieties,  Invent. Math.  159  (2005), no.
3, 607--656.
\item[{\bf [Ja1]}] U. Jannsen, On the $\ell$-adic cohomology
of varieties over number fields and its Galois cohomology, in {\it
Galois Groups over \bQ}, Y.Ihara, K. Ribet and J.-P. Serre ed.,
Math. Sci. Res. Inst. Publ. 16, Springer-Verlag, New York (1989),
315-360.
\item[{\bf [Ja2]}]  U. Jannsen, Mixed Motives and Algebraic K-Theory,
Lecture Notes in Mathematics, Volume 1400, Springer-Verlag 1990
\item[{\bf [Ka]}]  K. Kato, Logarithmic structures of Fontaine-Illusie, in
Algebraic Analysis, Geometry and Number Theory, Proceedings of the
JAMI Inaugural Conference, Jun-Ichi Igusa, editor, Special Issue of
the American Journal of Mathematics, Johns Hopkins Press, 1989
\item[{\bf [Ku1]}]  K. K{\"u}nnemann, Projective regular models for abelian
varieties, semi-stable reduction, and the height pairing, Duke Math.
J. 95 (1998) 161-212.
\item[{\bf [Ku2]}]  K. K{\"u}nnemann, Algebraic cycles on toric
fibrations over abelian varieties, Math. Zeitschrift 232 (1999)
427-435
\item[{\bf [La]}]  R. Langlands, Sur la mauvaise r{\'e}duction d'une
vari{\'e}t{\'e} de Shimura, in Journ{\'e}es de G{\'e}om{\'e}trie
Alg{\'e}brique de Rennes, Ast{\'e}risque 65 (1979)
\item[{\bf [Mok]}] A. Mokrane, La suite spectrale des poids
en cohomologie de Hyodo-Kato, Duke Math. J. 72 (1993) 301-337.
\item[{\bf [Mu1]}]  D. Mumford, An analytic construction of degenerating
curves over a complete local ring, Compositio Math. 24 (1972)
129-174.
\item[{\bf [Mu2]}]  D. Mumford, An analytic construction of degenerating
abelian varieties over complete local rings, Compositio Math. 24
(1972) 239-272.
\item[{\bf [Mus]}]  G.A. Mustafin, Non archimedean uniformization (in
Russian),
 Mat. Sb. 105 (147) (1978) 207-237, 287; English translation in Mathematics
of the USSR:  Sbornik 34 (1978) 187-214.
\item[{\bf [Na]}]  C. Nakayama, Degeneration of $l$-adic weight spectral
sequences, Amer. J. Math. 122 (2000), no. 4
\item[{\bf [P]}] B. Perrin-Riou, Repr{\'e}sentations $p$-adiques
ordinaires, in {\it P{\'e}riodes
 p-adiques, S{\'e}minaire de Bures, 1988} (J.-M. Fontaine ed.),
 Ast{\'e}risque 223 (1994), 185-208.
\item[{\bf [Rap]}]  M. Rapoport, On the bad reduction of Shimura varieties,
in Automorphic Forms, Shimura Varieties and $L$-Functions, L. Clozel
and J. S. Milne editors, Perspectives in Mathematics, Vol. 11,
Academic Press 1990
\item[{\bf [RZ1]}] M. Rapoport and Th. Zink, {\"U}ber die locale
Zetafunktion von Shimuravariet{\"a}ten, Monodromie-filtration und
verschwindende Zyklen in ungleicher Characteristic, Invent. Math. 68
(1982), 21-101.
\item[{\bf [RZ2]}] M. Rapoport and Th. Zink, Period spaces for
$p$-divisible groups, Annals of Mathematics Studies, Volume 141,
Princeton University Press, 1996
\item[{\bf [R1]}]  W. Raskind, The $p$-adic Tate conjecture for divisors on varieties
with totally degenerate reduction, in preparation 2004
\item[{\bf [R2]}]  W. Raskind, A generalized
Hodge-Tate conjecture for varieties with totally degenerate
reduction over $p$-adic fields, preprint 2004, to appear in the
Proceedings of the International Conference on Algebra and Number
Theory, Hyderabad, India
\item[{\bf [RX]}] W. Raskind and X. Xarles, On $p$-adic
intermediate jacobians, to appear in Trans. Amer. Math. Soc.
\item[{\bf [Ray]}] M. Raynaud, Vari{\'e}t{\'e}s ab{\'e}liennes et
g{\'e}om{\'e}trie rigide, Actes du ICM, Nice, 1970, Tome I, 473-477
\item[{\bf [ST]}]  J.-P. Serre and J. Tate, Good reduction of
abelian varieties, Annals of Math. 88 (1968) 492-517
\item[{\bf [St]}] J.H. Steenbrink, Limits of Hodge Structures,
Invent. Math. 31 (1976), 223-257.
\item[{\bf [Tsu]}] T. Tsuji, $p$-adic {\'e}tale cohomology and
crystalline cohomology in the semistable reduction case, Inventiones
Math 137 (1999) 233-411.
 \item[{\bf [Va]}]  Y. Varshavsky, $p$-adic uniformization of unitary
Shimura varieties, Publ. Math. I.H.E.S. 87 (1998) 57-119.
\item[{\bf [Z]}]  T. Zink, {\"U}ber die schlechte Reduktion einiger
Shimuramannig\-faltig\-kei\-ten, Comp. Math 45 (1981) 15-107.
\end{itemize}

\textit{Authors' addresses}
\vspace{.1in}\\
\textsc{Wayne Raskind:  Department of Mathematics\\University of
Southern California\\Los
Angeles, CA 90089-2532, USA}\\
\texttt{email:  raskind@math.usc.edu}\\

\textsc{Xavier Xarles:  Departament de Matem{\`a}tiques\\Universitat
Aut{\`o}noma de
Barcelona\\08193 Bellaterra, Barcelona, Spain}\\
\texttt{email:  xarles@mat.uab.es}

\end{document}